\nonstopmode

\documentclass[11pt]{article}

\usepackage[fleqn]{amsmath}

\usepackage{hyperref}
\hypersetup{
    colorlinks=true,
    linkcolor=blue,
    filecolor=magenta,      
    urlcolor=blue,
}

\usepackage{tcolorbox}   
\definecolor{mycolor}{rgb}{0.122, 0.435, 0.698}

\usepackage{amsmath,amssymb,amsthm}
\usepackage{graphics,epsfig,calc}

\usepackage{latexsym,epsfig,bm,amssymb}
\usepackage{xcolor}
\usepackage{amsthm,mathrsfs}

\usepackage{mathptmx}

\DeclareSymbolFont{AMSb}{U}{msb}m{n}
\DeclareSymbolFontAlphabet{\mathbb}{AMSb}

\newcommand{\beqn}{\begin{eqnarray}}
\newcommand{\eeqn}{\end{eqnarray}}
\newcommand{\be}{\begin{equation}}
\newcommand{\ee}{\end{equation}}

\newcommand{\bsp}{\begin{split}}
\newcommand{\esp}{\end{split}}

\newcommand{\ba}{\begin{array}}
\newcommand{\ea}{\end{array}}

\newcommand{\bpr}{\begin{proof}}
\newcommand{\epr}{\end{proof}}

\newcommand{\bH}{{\bf H}}

\newcommand{\bA}{{\bf A}}

\newcommand{\cD}{{\cal D}}

\newcommand{\bbd}{{\bf d}}

\newcommand{\bF}{{\bf F}}

\newcommand{\cH}{{\cal H}}

\newcommand{\bM}{{\bf M}}

\newcommand{\bN}{{\bf N}}

\newcommand{\cO}{{\cal O}}

\newcommand{\bS}{{\bf S}}
\newcommand{\cT}{{\cal T}}

\newcommand{\cU}{{\cal U}}

\newcommand{\aX}{{\mathbb X}}

\newcommand{\aM}{{\mathbb M}}
\newcommand{\aS}{{\mathbb S}}

\newcommand{\bX}{{\bf X}}

\newcommand{\bZ}{{\bf Z}}

\newcommand{\al}{\alpha}

\newcommand{\ci}{\cite}
\newcommand{\de}{\delta}
\newcommand{{\De}}{{\Delta}}

\newcommand{\ds}{\displaystyle}
\newcommand{\fr}{\frac}

\newcommand{\ga}{\gamma}

\newcommand{\la}{\label}
\newcommand{\lam}{\lambda}

\newcommand{  \om}{  \omega}
\newcommand{  \Om}{  \Omega}

\newcommand{\vp}{\varphi}

\newcommand{ \ov}{ \overline}

\newcommand{\pa}{\partial}

\newcommand{\trace}{{\rm tr\5}}

\newcommand{\we}{\wedge}

\newcommand{\si}{\sigma}
\newcommand{\ti}{\tilde}
\newcommand{\dist}{{\rm dist\5}}

\newcommand{\ve}{\varepsilon}

\newcommand\C{{\mathbb C}}
\newcommand\R{{\mathbb R}}

\newcommand{\bB}{{\bf B}}

\newcommand{\vka}{\varkappa}

\newcommand{\cm}{{\rm m}}

\newcommand{\5}{{\hspace{0.5mm}}}

\newcommand{\1}{{\hspace{0.1mm}}}

\newcommand{\const}{\mathop{\rm const}\nolimits}

\newcommand{\tr}{\mathop{\rm tr\,}\nolimits}

\newcommand{\rRe}{{\rm Re\5}}
\newcommand{\rIm}{{\rm Im\5}}

\newcommand{\Ker}{{\rm Ker\5}}

\newtheorem{theorem}{Theorem}[section]

\renewcommand{\thetheorem}{\arabic{section}.\arabic{theorem}}

\newtheorem{definition}[theorem]{Definition}

\newtheorem{lemma}[theorem]{Lemma}
\newtheorem{example}[theorem]{Example}

\newtheorem{remark}[theorem]{Remark}

\newtheorem{remarks}[theorem]{Remarks}
\newtheorem{cor}[theorem]{Corollary}
\newtheorem{proposition}[theorem]{Proposition}

\newcommand{\bd}{\begin{definition}}
 \newcommand{\ed}{\end{definition}}

\newcommand{\bt}{\begin{theorem}}
 \newcommand{\et}{\end{theorem}}
\newcommand{\bqt}{\begin{qtheorem}}
 \newcommand{\eqt}{\end{qtheorem}}

\newcommand{\bp}{\begin{proposition}}
 \newcommand{\ep}{\end{proposition}}

\newcommand{\bl}{\begin{lemma}}
 \newcommand{\el}{\end{lemma}}
\newcommand{\bc}{\begin{cor}}
 \newcommand{\ec}{\end{cor}}

\newcommand{\bex}{\begin{example}}
 \newcommand{\eex}{\end{example}}
 
\newcommand{\bexs}{\begin{examples}}
 \newcommand{\eexs}{\end{examples}}

\newcommand{\bexe}{\begin{exercice}}
 \newcommand{\eexe}{\end{exercice}}

\newcommand{\br}{\begin{remark}}
 \newcommand{\er}{\end{remark}}
 
\newcommand{\brs}{\begin{remarks}}
 \newcommand{\ers}{\end{remarks}}

\newcommand{\bce}{\begin{center}}
\newcommand{\ece}{\end{center}}

\begin{document}
\begin{center}

{\huge On single-frequency asymptotics for  
\bigskip

the Maxwell--Bloch equations: mixed states
}
\bigskip\smallskip

 {\large A.I. Komech$^1$ and E.A. Kopylova}\footnote{ 
 Supported by Austrian Science Fund (FWF) PAT 3476224}
 \\
{\it
Institute of Mathematics of
BOKU
University, Vienna, Austria\\
}
 alexander.komech@boku.ac.at,  elena.kopylova@boku.ac.at

\end{center}

\setcounter{page}{1}
\thispagestyle{empty}

\begin{abstract}
We consider damped driven
Maxwell--Bloch equations   which are  finite-dimensional 
approximation of the damped driven  Maxwell--Schr\"odinger equations. The  
 equations describe a single-mode Maxwell field
coupled to a two-level molecule. 
  Our main result is the construction of solutions with single-frequency asymptotics of the 
Maxwell field in the case of   quasiperiodic pumping.
The asymptotics hold for solutions with {\it harmonic} initial values
which are stationary states of averaged  equations in the interaction picture.

We
calculate all  harmonic states
and analyse  their stability.
The calculations rely on the Bloch--Feynman gyroscopic representation of 
von Neumann equation for the density matrix.
The asymptotics follow by application of
the  averaging  theory of the Bogolyubov type.
The key role 
 in the application of the averaging theory
 is played by  a special a priori estimate.

  \end{abstract}
  
  \smallskip
  

\tableofcontents

\setcounter{equation}{0}
\section{Introduction}
The Maxwell--Bloch equations (MBE) were introduced by Lamb \ci{L1964}
for the semiclassical 
description of the laser action \ci{AE1987,BM2011,H1984,SSL1978,SZ1997,S1986,S2012}.
The equations
are the Galerkin approximation of the Maxwell--Schr\"odinger  system
 \ci{BT2009,GNS1995,K2019phys,K2013,K2022,NW2007,S2006}
 Our main goal is construction of solutions with single-frequency asymptotics. The asymptotics seem to correspond to 
 the laser coherent radiation, which
remains a key mystery of laser action since its discovery around 1960.
The damped driven  MBE for mixed states read
\be\la{HMB2} 
\left\{\ba{rcl}
\dot A(t)\!\!&\!\!=\!\!&\!\!B(t),\quad
\dot B(t)=-\Om^2 A(t)\!-\!\ga  B(t)+cj(t)
\\\\
i\hbar\dot \rho(t)\!\!&\!\!=\!\!&\!\![H(t),\rho(t)]
 \ea\right|,\qquad j(t)={2}\vka\,\rIm\rho_{21}(t).
\ee
Here $A(t),B(t)\in\R$,  $\Om>0$ is the resonance frequency, 
$\ga>0$ is the dissipation coefficient,
$c$~is the speed of light, and $\hbar$ - the Planck constant.
Further,
$H(t)$ is the Hermitian matrix, and 
$\rho(t)=\rho$ is a
 nonnegative Hermitian  $2\times 2$-density matrix,  
\be\la{Hro}
H(t)=
\left(\!\!\!\ba{cc}\hbar\om_1&i a(t)\\\\-ia(t)&\hbar\om_2\ea\!\!\!\right);\qquad
\rho=
\left(\!\!\!\ba{cc}\rho_{11}&\ov\rho_{21}\\\\ \rho_{21}&\rho_{22}\ea\!\!\!\right),
\qquad \tr\rho=1,\quad\rho\ge 0,
\ee
where $\hbar\om_2>\hbar\om_1$ are the energy levels of active molecules, and
 the function $a(t)$ is given by
\be\la{JOm}
a(t):=\fr {\vka}c [A(t)+A^e(t)],\qquad \vka=p\,{{\om}},\quad 
{\om}=\om_2-\om_1>0,
\ee
where $p\in\R$ is proportional to the molecular dipole moment.
The pumping  $A^e(t)$ we suppose to be a quasiperiodic function:
\be\la{qpp}
A^e(t)=\rRe[\bA^e e^{-i\Om t}]+\rRe\sum_1^N \bA^e_k e^{-i\Om_k t},
\qquad{\rm where}\qquad \bA^e,\bA^e_k\in\C,\quad \Om_k\in\R\quad{\rm and}\quad 
\Om_k\ne\Om.
\ee
For solutions to the MBE,
the   conservation 
 $\trace\rho(t)=\const$
 holds 
  since the trace of commutators vanish. We will consider
solutions with $\const=1$; see (\ref{Hro}).
We use the Heaviside--Lorentz units and
recall the introduction of the equations in Appendix \ref{sMB}.

  In Appendix \ref{swp}, we
  prove
  the a priori bound  for solutions $X(t)=(A(t),B(t),\rho(t))$ 
  to  (\ref{HMB2}) with $|p|/\ga=r$
 satisfy the  following a priori bounds: 
   \be\la{apri}
|X(t)|\le D_r(|X(0)|),\qquad t\ge 0,
\ee
which is proved in  Appendix \ref{swp}.
The bound  implies
  the well-posedness of the MBE.
  \smallskip

 Our main goal is  asymptotics of the Maxwell amplitudes $A(t)$ and $B(t)$
    for solutions $X(t)$ to the  MBE
  as
  $|p|,\ga\to 0$.
  Define   complex Maxwell amplitudes by $M(t)=A(t)+iB(t)/\Om$.
Then the first two equations of the MBE are equivalent to
\be\la{Max}
\dot M(t)=-i \Om M(t)-i\ga  M_2(t)+2icp\om\,\rIm\rho_{21}
%
/\Om,\quad{\rm where}\quad M_2(t)=\rIm M(t).
\ee

Note that the parameters $|p|,\ga$ are very small for many types of lasers, see Appendix \ref{aC}. 
For $p=\ga=0$,  all solutions to (\ref{Max})
are single-frequency: $M(t)=e^{-i\Om t}M(0)$.
For small $p$ and $\ga$, equation (\ref{Max}) implies that 
$
M(t)= e^{-i\Om t}M(0) +\ds\int_0^t R(s)ds,
$
 where $\sup\limits_{s\ge 0} |R(s)|=\cO(|p|+\ga )$  by  (\ref{apri}). Hence,  $M(t)= e^{-i\Om t}M(0) +\cO(|p|t)$
for $p/\ga=r$ with arbitrary  $r\ne 0$.
In particular, for solutions with any fixed initial states $M(0)$,
 \be\la{intriad0}
 \max_{t\in [0,|p|^{-1/2}]}|M(t)-e^{-i\Om t}M(0)|=\cO({|p|^{1/2}}),
 \qquad { p\to 0},\quad{p}/\ga=r.\quad\qquad\qquad\qquad\qquad
   \ee 
   Our main results show that in the resonance case, when
  $\Om=\om$,
the time scale $p^{-1/2}$ in this asymptotics can be extended 
to $|p|^{-1}$ and even more in the case of special  {\it harmonic states}
of the MBE.
\smallskip  

 We will represent the density matrix via the Pauli matrices as
\be\la{roS}
\rho(t)=\fr12[E+S(t)\cdot \si],\qquad {\rm where} 
\quad S(t)\in\bB=\{S\in\R^3: |S|\le 1\},\quad \si=(\si_1,\si_2,\si_3),
\ee
and accordingly, solutions $X(t)$ to MBE are represented by 
$(M(t),S(t))\in \aX=\C\times\bB$.
     The harmonic states are defined via the
 {\it  interaction picture} of the MBE
 and 
 the corresponding {\it averaged equation} with the structure
 \be\la{aver}
 \dot \bX(t)=p\bF_r(\bX(t)),\qquad t\ge 0;\qquad r=p/\ga.
 \ee
 We define 
  the harmonic states $\bX\in\aX$ of the MBE 
   as stationary solutions of (\ref{aver}).
  
  We
  calculate all  harmonic states $\bX^r=(\bM^r,\bS^r)\in\aX$ of the MBE with $p/\ga=r$ and show that
   the set  is a  union $Z^r=Z^r_1\cup Z^r_2$ of two smooth 1D manifolds.
  We linearise the   dynamics (\ref{aver}) at the harmonic states and 
  calculate the spectra of the linearisations.
  The  calculations show that
  for $\Om=\om$ and
$c|r|>|\bA^e|$, there exist a 
 nonempty submanifold $Z^r_+\subset Z^r$ which is attractive under the averaged dynamics (\ref{aver}).
  Denote the matrix
   \be\la{VOm}
 \qquad \qquad \qquad \qquad  V_\Om=\left(\ba{ccc}
0&\Om&0\\
-\Om&0&0\\
0&0&0\ea
\right).
   \ee
   {
   Our main results are the following asymptotics
   for solutions $X(t)=(M(t),S(t))$ to equations (\ref{HMB2}) in the representation (\ref{roS}),
   with a
fixed quotient $p/\ga=r$. Everywhere below we  consider the case $p>0$ only  since the extension
to $p<0$ is obvious.

\bt\la{t1} 
Let $\Om=\om$ and the pumping be quasiperiodic. Then for any $r>0$,
the following asymptotics hold.
\smallskip\\
i) Let the  initial state $X(0)=(\bM,\bS)\in Z^r$. Then 
the 
corresponding 
solutions 
admit the following {\it adiabatic asymptotics}:
 \be\la{intriad}
 \max_{t\in [0,p^{-1}]}\Big[|M(t)-e^{-i\Om t} \bM
 +|S(t)-e^{V_\Om t} \bS|\Big]
 =\cO(p^{1/2}),
 \qquad p\to 0.\qquad\qquad\qquad
 \quad
   \ee      
ii)
 Let 
$cr>|\bA^e|$,  and  $D^r_ d$ denote a suitable  subset of the 
tubular $ d$-neighborhood  of the stable submanifold $Z^r_+$
with sufficiently small $ d>0$. Then 
for initial states $X(0)\in D^r_ d$, 
the  corresponding solutions 
admit the following asymptotics uniformly 
in initial values $X(0)\in D^r_ d$:
 \be\la{intrias}
 \max_{t\in [0,p^{-1}]}\Big[ |M(t)-e^{-i\Om t}\bM_*|+|S(t)-e^{V_\Om t} \bS_*(t)|\Big]
 =\cO[{p^{1/2}}+ d],
 \qquad p\to 0,\qquad
   \ee 
   where $(\bM_*,\bS_*(t))\in Z^r_+$ and 
       $\bM_*=-\bA^e$ does not depend on $r$ and on $X(0)\in D^r_ d$. 
       The Lebesgue measure 
   \be\la{slo}
   |D^r_d|\sim  d^4,\qquad  d\to 0.
\ee
  \et
  }
  \br\rm
   Without 
   the constraint $p/\ga=r$, the
       a priori bounds (\ref{apri}) and the asymptotics (\ref{intriad})--(\ref{intrias}) do not hold.
  For example,  the bounds must be different for small and big $|p|$ for a given dissipation $\ga>0$.
  The asymptotics 
   do not hold without this restriction since the limiting amplitudes $\bM^r$ depend on $r$.

   \er

 The asymptotics (\ref{intriad}) for $S(t)$ reads
 \be\la{asS}
 S(t)= e^{V_\Om t} \bS+\cO(p^{1/2}),\qquad t\in[0,p^{-1}]
  \ee
  which is approximately the precession of $S(t)$ about axis $S_3$.
  
  \smallskip
  
Let us comment on
 our approach. 
 The asymptotics (\ref{intriad}) mean that
 the MBE admits solutions  
$M(t)=e^{-i\Om t}\aM^r(t)$, $S(t)=e^{V_\Om t}\aS^r(t)$ 
with {\it slowly varying} enveloping amplitudes $\aM^r(t)\in\C$ and $\aS^r(t)\in\R^3$
for small $p$ and $\ga$ with $p/\ga=r$.
The amplitudes are solutions to  the corresponding dynamical system
which is the interaction picture (or ``rotating frame representation") of the MBE.
 The slow variation of the amplitudes for small $p$ and $\ga$ is equivalent to the fact  that
the initial state $(\aM(0),\aS(0))=(M(0),S(0))$
is a harmonic state, i.e.,
a stationary solution of (\ref{aver}) with $p/\ga=r$.

 The calculation  of harmonic states 
 and of spectra of the linearisations
 rely on  
 the Bloch--Feynman representation  \ci{B1946,FVH1957} 
  of the von Neumann equation via the vector
   $S(t)$ from (\ref{roS}).
   We prove that the harmonic states $(\bM,\bS)$ with $\bM\ne 0$
  exist only in the resonance case $\Om=\om$, 
and  only
 if $\bA^e\ne 0$.
 In particular, for the single-frequency pumping
 $A^e(t)=\rRe [\bA_p e^{-i\om_p t}]$,
 the harmonic states with $\bM\ne 0$ exist only in the case $\bA_p\ne 0$
 and triple resonance
 \be\la{trip} 
  \Om=\om=\om_p.
  \ee

 The asymptotics (\ref{intriad}) follows for solutions with
 every initial state $(M(0),S(0))\in Z^r$ 
 by the averaging theory of Bogolyubov type  \ci[Theorem 4.3.6]{SVM2007}.
To prove asymptotics (\ref{intrias}), we combine 
the theorem with
the stability of the branch  $Z^r_+$.
 Our calculations show that for all harmonic states $(\bM,\bS)\in Z^r_+$,
the component $\bM=-\bA^e$. This is why
the limiting amplitude $\bM_*$ in (\ref{intrias}) does not depend on
$r>0$ and on
 $X(0)\in D^r_d$. The key role 
 in the application of the averaging theory
 is played by the special a priori estimate (\ref{apri}).

\br\rm
 Our approach relies on the averaging theory which
neglects oscillating terms. So, it gives a justification of the ``rotating wave approximation",  which is widely used in  Quantum
Optics \ci{AE1987,H1984,SSL1978,SZ1997,S1986,S2012}.
The asymptotics (\ref{intriad}), (\ref{intrias}) 
specify the  time scale and  the  error of such approximations.
\er

Let us comment on related  results.
The problem of existence of time-periodic solutions 
to the MBE 
has been
 discussed since 1960s. 
 The first results in this direction were obtained recently in
  \ci{CLV2019} and \ci{WWG2018}
 for various versions
 of the MBE. 
   In  \ci{WWG2018},
  the  N-th order time-periodic
  solutions 
  were constructed by perturbation techniques.
    For the phenomenological model \ci{A1985,AB1965},
  time-periodic solutions
  were constructed in \ci{CLV2019} in the absence of
   time-periodic pumping  for small interaction constants.
   The solutions are obtained as the result of a bifurcation
   relying on  homotopy invariance of the degree  \ci{CLN2017} and developing the
   averaging arguments  \ci{CL2018}. 
     The period is determined by bifurcation.
     
In \ci{K2024}, we have established the existence of 
solutions with $T$-periodic Maxwell amplitude
for any $T$-periodic pumping 
without smallness conditions.     
     
     In \ci{KK2025sin}, we have constructed solutions with asymptotics 
of type (\ref{intriad}), (\ref{intrias}) for the MBE  with {\it pure states},  
using the reduction by the symmetry gauge group. 
In this case,
for any $r>0$, the set of stationary states is discrete 
and consists of one or two points. { The construction of the asymptotics for  
 system (\ref{HMB2}) in the present paper, 
 required completely new technique based on the Bloch--Feynman
 representation.
 }
 
Up to our knowledge,
the single-frequency asymptotics  for the MBE with mixed states were
not constructed till now.

\smallskip

Let us comment on our exposition.
In Section \ref{sA}, we construct the representation of the Bloch--Feynman
type for the von Neumann equation from MBE.
The dynamics on the 
interaction picture is calculated in Section \ref{sav}, its averaging in Section \ref{saver},
and  all stationary states of the averaged equations are calculated in
Section \ref{sst}.
In Section \ref{sstab} we analyse the stability of the stationary states, and 
in Section \ref{ssin} we prove the single-frequency asymptotics 
(\ref{intriad})--(\ref{intrias}).
In Appendix \ref{swp} we establish the bound (\ref{apri}).
In
Appendix \ref{sMB}
 we comment on  the introduction of the MBE with mixed states, and
 in Appendix \ref{aC} we discuss possible treatment of the laser threshold
 and laser amplification relying on our results.

{\bf Acknowledgements.} The authors thank S. Kuksin, M.I. Petelin, A. Shnirelman and H. Spohn
 for longterm fruitful discussions, and the  Institute of Mathematics of BOKU University
for the support and hospitality. The research is supported by
Austrian Science Fund (FWF) PAT 3476224.

\setcounter{equation}{0}

\section{The Bloch-Feynman ``gyroscopic" representation}\la{sA}
Here we represent von Neumann equation from MBE
as a ``gyroscopic equation" (\ref{ome})
applying the Bloch, Feynman \& al approach \ci{BM2011,B1946,FVH1957,GK2004}.
Hermitian density matrices $\rho$ with $\tr\rho=1$
admit the representation
\be\la{rhob}
\rho=\fr12\left(\!\!\ba{cc}1+S_3&S_1-iS_2\\S_1+iS_2&1-S_3\ea\!\!\right),\qquad 
(S_1,S_2,S_3)\in\R^3.
\ee
Hence, using the Pauli matrices and the nonnegativity $\det\rho=1-S_1^2-S_2^2-S_3^2\ge 0$, we obtain (\ref{roS}):
\be\la{BS}
\rho=\fr12\left(\!\!\ba{cc}1&0\\0&1\ea\!\!\right)
+
\fr{S_1}2\left(\!\!\ba{cc}0&1\\1&0\ea\!\!\right)
+\fr{S_2}2\left(\!\!\ba{cc}0&-i\\i&0\ea\!\!\right)
+\fr{S_3}2\left(\!\!\ba{cc}1&0\\0&-1\ea\!\!\right)=\fr12[E+
S\cdot \si].
\ee
In particular, we can expand the ``Hamiltonian operator"
(\ref{Hro}) as
\be\la{Ham2}
H(t)=\fr\hbar2(\om_1+\om_2)E+\fr\hbar2(\om_1-\om_2)\si_3-a(t)\si_2
=\fr\hbar2(\om_1+\om_2)E+\fr\hbar2(\om_1-\om_2)e_3\cdot \si-a(t)e_2\cdot\si,
\ee
where $e_3=(0,0,1)$ and $e_2=(0,1,0)$.
Now the von Neumann equation in MBE for $\rho(t)=\fr12[E+S(t)\cdot\si]$  reads as
\be\la{vN}
i\hbar \dot S(t)\cdot\si=2i\fr\hbar2(\om_1-\om_2)[e_3\we S(t)]\cdot\si-a(t)[e_2\we S(t)]\cdot\si=2i\Big[\big(\fr\hbar2(\om_1-\om_2)e_3-a(t)e_2\big)\we S(t)\Big]\cdot\si,
\ee
where we have used the formula $[a\cdot\si,b\cdot\si]=2i(a\we b)\cdot\si$ for $a,b\in\R^3$.
This implies
the  ``gyroscopic equation" of type  \ci[(4)]{FVH1957}:
\be\la{ome}
 \dot S(t)=\theta(t)\we S(t),\qquad \theta(t)=- \om \,e_3-2\fr{a(t)}\hbar e_2,\quad\om=\om_2-\om_1.
\ee
The equation can be also written  as
\be\la{vNma}
\dot S(t)=\Theta(t) S(t),\qquad \Theta(t)=
\begin{pmatrix}
0&-\theta_3(t)&\theta_2(t)\\
\theta_3(t)&0&-\theta_1(t)\\
-\theta_2(t)&\theta_1(t)&0\end{pmatrix}=\begin{pmatrix}
0&\om&-2 a(t)/\hbar\\
-\om&0&0\\
2 a(t)/\hbar&0&0\end{pmatrix}\in so(3),
 \ee
 which corresponds to  the  cross-product representation
 of the Lie algebra $so(3)$.
 In particular,  the following conservation law holds,
 \be\la{Sc}
 |S(t)|=\const,\qquad t\ge 0.
 \ee
%
The current $j(t)=\vka S_2(t)$, so
the MBE (\ref{HMB2}) reduces to (cf. (\ref{Max}))
\be\la{HMB2M} 
\left\{\ba{rcl}
\dot M(t)\!\!&\!\!=\!\!&\!\!-i \Om M(t) -i\ga  M_2(t)+ic   \vka S_2(t)/\Om
\\\\
\dot S(t)\!\!&\!\!=\!\!&\!\! {\Theta(t)} S(t)
 \ea\right|,\qquad M(t)=
A(t)+iB(t)/\Om.
\ee

  \setcounter{equation}{0}
  \section{Dynamics in the interaction picture}\la{sav}
   By (\ref{vNma}), for small $p,\ga>0$,
the system (\ref{HMB2M}) is a small perturbation of the unperturbed one,
\be\la{systu} 
\dot M( t)=-i \Om M( t),
\qquad
\dot S( t)=V_\om S( t),
\ee
where $V_\om$  is the matrix 
(\ref{VOm})  with $\Om=\om$ which coincides with 
(\ref{vNma}) in the case $a(t)=0$.
Solutions to this system are given by 
\be\la{soli}
 M( t)=e^{-i  \Om t}\aM, \qquad S( t)= e^{V_\om  t}\aS,\qquad{\rm where} \qquad \aM\in\C,\quad \aS\in \bB.
 \ee
Our goal is  construction of similar solutions
to the perturbed system (\ref{HMB2M}),
 \be\la{intr}
   M(  t)=e^{-i  \Om t}\aM(  t), \qquad S(  t)=e^{V_\om  t}\aS(  t)
   \ee
with slowly varying enveloping amplitudes: for a wide interval of time $[0,T(p)]$
\be\la{intr2}
 \sup_{ t\in [0,T(p)]}  \big[|\aM(  t)- \aM(0)|+| \aS(  t)-\aS(0)|\big]\to 0,\qquad p\to 0, \quad p/\ga=r\ne 0.
   \ee
Substituting (\ref{intr}) into (\ref{HMB2M}), we obtain  the dynamical equations
for the enveloping amplitudes (the ``interaction picture"):
 \be\la{systim}
  \left\{\ba{rl}
  \dot \aM( t)&=-ie^{i \Om t}\big[\ga\rIm( e^{-i\Om t}\aM( t)) -\ti\vka ( e^{V_\om  t}\aS(  t))_2\big]\\\\
  \dot \aS(  t)&=e^{-V_\om  t}({\Theta( t)}-V_\om)e^{V_\om  t}\aS( t)
  \ea\right|,\qquad \ti\vka=\fr{c\vka}\Om.
  \ee
   The equations are called as the interaction picture of (\ref{HMB2M}).
  By (\ref{ome}),  (\ref{vNma}), and (\ref{JOm}),  
  \be\la{enu1}
  \Theta( t)-V_\om=
 \fr{2a( t)}\hbar 
  \left(\ba{ccc}
0&0&-1\\
0&0&0\\
1&0&0\ea
\right),\qquad a( t):=\fr {\vka}c [A( t)+A^e( t)],\qquad \vka=p\,{{\om}}.
 \ee
  Note that $e^{V_\om  t}$ is he dynamical group of the gyroscopic equation (\ref{ome}) with 
  the angular velocity 
  $\theta( t)=-\om e_3$. Hence, $e^{V_\om  t}$ is the rotation about $e_3$ with
   the angular velocity
  $\om$:
  \be\la{enu2}
  e^{V_\om  t}=\left(\ba{ccc}
\,\,\,\,\,\,\,\,\cos\om t &\sin\om t&0\\
-\sin\om t&\cos\om t&0\\
0&0&1\ea
\right) 
  \ee
    Substituting (\ref{enu1}) and  (\ref{enu2}) into   (\ref{systim}),  we obtain
 $$
  \left\{\ba{rl}
     \dot \aM( t)=\!\!&\!\!- p ie^{i\Om  t}\!\big[\ga_1
   (\aM_2( t)\cos\Om    t-\aM_1( t)\sin \Om   t)
 -\vka_1 
 ( -\aS_1( t)\sin\om t+\aS_2( t)\cos\om t)
  \big],\\\\
  \dot \aS(  t)=\!\!&\!\!- pb [A( t)\!+\!A^e( t)] \big[\aS_3( t)( e_1\cos\om t+ e_2\sin\om t)
  -e_3(\aS_1( t)\cos\om t+\aS_2( t)\sin\om t) \big],
 \ea\right.
$$
where 
 \be\la{P1}
\ga_1=\fr \ga {p}=\fr1{ r},\qquad \vka_1=\fr{\ti\vka}p=\fr{c\,{\om}}{\Om},\qquad  b=\ds\fr{2\om}{c\hbar}.   
  \ee
  The  equations can be written as
  \be\la{perN}
  \dot \aM( t)=p \,f_r(\aM( t),\aS( t), t),\qquad \dot \aS( t)=p \,g_r(\aM( t),\aS( t), t),
  \ee  
 where the functions $f$ and $g$ are given by 
 \be\la{fg}
 \!\!\!\!
 \left\{\!\!\!\!\!\ba{rl}
 f_r(\aM,\aS, t)\!\!=\!\!\!&\!\!\! \!\!\big[\!-\!i\cos \Om  t\!+\!\sin\Om   t\big]
  \big[\ga_1(\aM_2\cos  \Om  t\!-\!\aM_1\sin  \Om  t)\! -\!\vka_1(\! -\aS_1\sin\om t\!+\!\aS_2\cos\om t)\big]
 \\\\
g_r(\aM,\aS, t)\!\!=\!\!\!&\!\!\!\!\!-b\big[\aM_1\cos\Om t\!+\!\aM_2\sin\Om t\!+\!A^e( t)\big] 
\big[\aS_3(e_1\cos\om t  \!+\!e_2 \sin\om t )\!-\!e_3(\aS_1\cos\om t\!+\!\aS_2\sin\om t) \!\big]
\ea\!\!\!\!\right|\!.\!\!
  \ee
  
  \br\rm
  It is important that 
  the coefficients $\ga_1,\vka_1,b$ depend only on 
  $\Om,\om, r$.
    Hence,
    for any fixed $\Om,\om>0$ and $r\ne 0$,
     the  asymptotics of solutions to systems(\ref{perN})
  as $p\to 0$ and $p/\ga=r$
   can be calculated by methods of the averaging  theory  \ci{B1961,SVM2007}.
  \er

  \setcounter{equation}{0}
  \section{The averaging}    \la{saver}
  
    The  
  averaged equations (\ref{perN}) read
    \be\la{perNa}
  \dot \bM( t)
  =p\,\ov f_r(\bM( t),\bS( t)),\qquad \dot \bS( t)=p\,\ov g_r(\bM( t),\bS( t)),
  \ee
 where 
  \be\la{fga}
  \ov f_r(\bM,\bS)=
  \langle
  f_r(\bM,\bS,\cdot)
  \rangle=
  \lim_{T\to\infty}\fr1T \int_0^T f_r(\bM,\bS, t)d t,\qquad
\ov g_r(\bM,\bS)=
  \langle
  g_r(\bM,\bS,\cdot)
  \rangle.
  \ee
  Let us calculate the averages (\ref{fga}). The results differ drastically 
for the resonance case $\Om=\om$ and non-resonance $\Om\ne \om$.
In notation (\ref{qpp}),   $\bA^e= \bA^e_1+i\bA^e_2$, where
\be\la{Apav}
 \bA^e
=2\langle A^e( t)e^{i \Om  t} \rangle,\qquad
  \bA^e_1
=2\langle A^e( t)\cos\Om   t \rangle,\qquad
  \bA^e_2=2\langle A^e( t)\sin\Om   t \rangle.
  \ee

\noindent
{\bf Resonance case $\Om=\om$.} Using the expressions (\ref{fg}), we obtain the averaged 
vector field
\be\la{avrnr}
 \ov f_r(\bM,\bS)=-\fr i2\![\ga_1\bM_2-\vka_1\bS_2\big]-  \fr12\big[\ga_1\bM_1-\vka_1\bS_1]
 = - \fr{\ga_1}2[\bM_1+i\bM_2]+ \fr{\vka_1}2[\bS_1+i\bS_2].
  \ee
\beqn\la{avr}
\ov g_r(\bM,\bS)&=&
-\fr{b}2\big[\bM_1[
  -e_3\bS_1 +\bS_3e_1]
+\bM_2\![
  -e_3\bS_2 +\bS_3e_2]
-\!e_3(\bS_1\bA^e_1+\bS_2\bA^e_2) +\!\bS_3(e_1\bA^e_1  \!+\!e_2 \bA^e_2)
\big]
 \nonumber \\
 \nonumber \\
  &=&
-\fr{b}2\big[
\bS_3
[e_1(\bM_1+\bA^e_1)+e_2(\bM_2+\bA^e_2)]
-e_3[\bS_1(\bM_1+\bA^e_1)+\bS_2(\bM_2+\bA^e_2)]\big].
\eeqn
\noindent      
{\bf Non-resonance case $\Om\ne \om$.} In this case, the calculations (\ref{avrnr}) simplify to
\be\la{avrnr2}
\ov f_r(\bM,\bS)=
 \ds-\fr{\ga_1}2 \bM,
 \ee
 which
implies 
nonexistence of stationary states with nonzero Maxwell field 
for the system (\ref{perNa}),
and also
the decay $\bM( t)=\ds\bM(0)e^{-\fr{\ga} 2 t}$.
  In contrast,
the resonance equation (\ref{avrnr})  includes the interaction term 
which can prevent the decay of
 the Maxwell amplitude, that is expected physically.
 
   \noindent{\bf The averaged system.}
By (\ref{avrnr}) and (\ref{avr}), the  averaged equations (\ref{perNa})
(or (\ref{aver}))  
in the resonance case
read as
\be\la{aveq}
\!\!\!\!\!\!\left\{\ba{rcl}
\dot\bM( t)&=&-\frac{p}2 \big[\ga_1[\bM_1+i\bM_2]-\vka_1[\bS_1+i\bS_2]\big]
\\\\
\dot\bS( t)&=&-\fr{pb}2\big[
\bS_3[e_1(\bM_1+\bA^e_1)+e_2(\bM_2+\bA^e_2)]
-e_3[\bS_1(\bM_1+\bA^e_1)+\bS_2(\bM_2+\bA^e_2)]
\big]
\ea\right|.
\ee
\setcounter{equation}{0}
\section{Stationary states for the  averaged dynamics (harmonic states)}\la{sst}
In this section, we calculate  all {\it harmonic  states} $(\bM,\bS)$,
i.e., stationary states
for the averaged equations (\ref{aveq}),
 in the resonance case $\Om=\om$. The states satisfy
\be\la{perNa1}
\left\{\ba{rl}
0=&
 \ds \ga_1[\bM_1+i\bM_2]- \vka_1[\bS_1+i\bS_2]
\\\\
0=&
\bS_3[e_1(\bM_1+\bA^e_1)+e_2(\bM_2+\bA^e_2)]-e_3[\bS_1(\bM_1+\bA^e_1)+\bS_2(\bM_2+\bA^e_2)]
\ea\right|.
\ee
It is important that the stationary equations depend on $r$ but {\it do not depend} on $p$.
The first equation of (\ref{perNa1}) gives
\be\la{al}
\bS_1=\al_r \bM_1,\quad\bS_2=\al_r \bM_2, \quad{\rm where}\quad \al_r=\fr{\ga_1}{\vka_1}=\fr1{cr}>0.
\ee
The
second equation of (\ref{perNa1}) is equivalent to the system
\be\la{sec}
\left\{\ba{rcl}
0&=&\bS_3(\bM_1+\bA^e_1)
\\
0&=&\bS_3(\bM_2+\bA^e_2)
\\
0&=&\bS_1(\bM_1+\bA^e_1)+\bS_2(\bM_2+\bA^e_2)
\ea\right|.
\ee
The last equation together with (\ref{al}) give $\bM_1(\bM_1+\bA^e_1)+\bM_2(\bM_2+\bA^e_2)=0$ which defines the circle \linebreak
$S:|\bM|^2+\bM\cdot\bA^e=0$, or equivalently, 
$|\bM+\bA^e/2|^2=|\bA^e|^2/4$, so
\be\la{cir}
\bM=\bM(\theta)=-\fr{\bA^e}2+\fr{|\bA^e|}2e^{i\theta},\qquad \theta\in[0,2\pi].
\ee
Recall that $|\bS|\le 1$ by (\ref{roS}).
Hence,
the set of all stationary states
with a fixed $r=p/\ga$
 is the union of two real  1D manifolds  $Z^r=Z^r_1\cup Z^r_2$, where 
\be\la{Sr1}
\!\!\!
Z^r_1=\{(\bM,\bS)\in \C\times \bB: \,\,\,
\bM\in S,
\,\,\,\,\bS_1=\al_r \bM_1,\,\,\,\,\bS_2=\al_r \bM_2,\,\,\,\,\bS_3=0\}.
\ee
The second   summand $ Z^r_2$  is 
\be\la{Sr2}
Z^r_2\!=\!\{(\bM,\bS)\in \C\times \bB: \,\,\,
\bM\!=\!-\bA^e,\,\,\,\bS_1\!=\!\al_r \bM_1,\,\,\,\bS_2\!=\!\al_r \bM_2\}.
\ee

\br\rm
i) $Z^r_1\ne\emptyset$ for all $r\ne 0$ since $\bM(0)=0$.
\smallskip\\
ii) $Z^r_2\ne\emptyset$ iff $|\al_r\bA^e|\le 1$, i.e. $cr\ge |\bA^e|$.
\smallskip\\
iii) The intersection of the manifolds is at most one point:
\be\la{S12}
Z^r_1\cap Z^r_2=\left\{\ba{ll}(-\bA^e,-\al_r\bA^e,0),&\al_r|\bA^e|\le 1\\ 
\emptyset,&\al_r|\bA^e|> 1\ea\right|.
\ee
\er
\br
\la{runif}
\rm
Let $(\bM(t),\bS(t))$ be a solution to 
the averaged
equations (\ref{aveq}). Then  $(e^{i\theta}\bM(t),e^{V_1\theta}
\bS(t))$ with $\theta\in\R$ is the  solution
to the same equations with $\bA^e$ replaced by $e^{i\theta}\bA^e$
and rotated vectors $e_1$ and $e_2$.
The same correspondence holds for solutions to
 stationary equations (\ref{perNa1})  that
is  particularly obvious  for the states (\ref{Sr2}).

\er

  \setcounter{equation}{0}
  \section{Spectra of linearised equations at the harmonic states}\la{sstab}
  Here, we calculate i) linearisations of the averaged system (\ref{aveq}) at the harmonic states and ii) their spectra. 
     By  (\ref{avrnr}) and (\ref{avr}), we obtain, omitting the index $r$:
$$
\left\{\ba{rl}
\ov f_i(\bM,\bS)=&- \frac 12\big[\ga_1\bM_i-\vka_1\bS_i\big], \quad i=1,2,
 \\\\
 \ov g_k(\bM,\bS) = & -\fr{b}2\big[\bS_3(\bM_k+\bA^e_k)\big],\quad k=1,2,
\\\\
 \ov g_3(\bM,\bS) = & \fr{b}2\big[\bS_1(\bM_1+\bA^e_1)+\bS_2(\bM_2+\bA^e_2)]
\ea\right.
$$
Differentiating, we get for $i,j,k,l=1,2$ 
$$
\left\{\!\!\!\!\ba{rclrclrcl}
\ds\frac{\pa \ov f_i}{\pa \bM_j}&=&-\frac{\ga_1}2\delta_{ij}, &\ds
\frac{\pa \ov f_i}{\pa  \bS_l}&=&\frac{\vka_1}2\delta_{il}, & 
\ds\frac{\pa \ov f_i}{\pa  \bS_3}&=&0
\\\\
\ds\frac{\pa \ov g_k}{\pa  \bM_j}&=&-\fr {b}{2}\bS_3\delta_{kj},&
\ds\frac{\pa \ov g_k}{\pa  \bS_l}&=&0, &\ds\frac{\pa \ov g_k}{\pa  \bS_3}&=&-\!\fr {b}{2}(\bM_k+\bA^e_k)
\\\\
\ds\frac{\pa \ov g_3}{\pa  \bM_j}&=&\fr {b}{2}\bS_j,&
\ds\frac{\pa \ov g_3}{\pa  \bS_l}&=&\fr {b}{2}(\bM_l+\bA^e_l), & \ds\frac{\pa \ov g_3}{\pa  \bS_3}&=&0
\ea\right.
$$
Thus, we obtain the Jacobian 
$$
J=J(\bM,\bS)=\frac 12\begin{pmatrix}
-\ga_1 &0&  \vka_1 & 0 & 0
\\ 
0& -\ga_1 &0&   \vka_1 &0 
\\
-b \bS_3 & 0 & 0& 0& -b(\bM_1+\bA^e_1)
\\
0&-b \bS_3  & 0& 0& -b(\bM_2+\bA^e_2)
\\
b\bS_1 &b\bS_2& b(\bM_1+\bA^e_1)  &b(\bM_2+\bA^e_2)& 0
\end{pmatrix}.
$$
Let us consider the cases $(\bM,\bS)\in Z^r_1$ and $(\bM,\bS)\in Z^r_2$ separately.
\smallskip\\
I. For $(\bM,\bS)\in Z^r_1$, the Jacobian simplifies to
$$
J=\fr12\begin{pmatrix}
-{\ga_1}& 0 & {\vka_1}&0&0
\\
0&-{\ga_1} & 0&{\vka_1}&0
\\
0&0&0&0&-b\bB_1
\\
0&0&0&0&-b\bB_2
\\
b\bS_1& b\bS_2&b\bB_1&b\bB_2&0
\end{pmatrix},
$$
where we denote $\bB_k=\bM_k+\bA^e_k$.  Hence,
$$
\det(2J-\lam)=
\begin{bmatrix}
-\ga_1-\lam & 0 & \vka_1&0&0
\\
0&-\ga_1-\lam & 0&\vka_1&0
\\
0&0&-\lam&0&-b\bB_1
\\
0&0&0&-\lam &-b\bB_2
\\
b\bS_1& b\bS_2&b\bB_1&b\bB_2&-\lam
\end{bmatrix}
$$
$$
=-(\ga_1+\lam) \begin{bmatrix}
-\ga_1-\lam & 0&\vka_1&0
\\
0&-\lam&0&-b\bB_1
\\
0&0&-\lam &-b\bB_2
\\
b\bS_2&b\bB_1&b\bB_2&-\lam
\end{bmatrix}
+\vka_1\begin{bmatrix}
0&-\ga_1-\lam & \vka_1&0
\\
0&0&0&-b\bB_1
\\
0&0&-\lam &-b\bB_2
\\
b\bS_1& b\bS_2&b\bB_2&-\lam
\end{bmatrix}
$$
$$
=(\ga_1+\lam)^2\begin{bmatrix}
-\lam&0&-b\bB_1
\\
0&-\lam &-b\bB_2
\\
b\bB_1&b\bB_2&-\lam
\end{bmatrix}
-\vka_1(\ga_1+\lam)\!\begin{bmatrix}
0&-\lam&-b\bB_1
\\
0&0&-b\bB_2
\\
b\bS_2&b\bB_1&-\lam
\end{bmatrix}
\!+\!\vka_1b\bS_1\!\begin{bmatrix}
\ga_1+\lam & \vka_1&0
\\
0&0&-b\bB_1
\\
0&-\lam &-b\bB_2
\end{bmatrix}
 $$
 $$
 =(\ga_1+\lam)^2\big(-\lam^3-\lam b^2\bB_1^2- \lam b^2\bB_2^2\big)
 -\vka_1(\ga_1+\lam)\lam b^2\bS_2\bB_2-\vka_1(\ga_1+\lam)\lam b^2\bS_1\bB_1
$$
$$
=-(\ga_1+\lam)\lam \big[(\ga_1+\lam)\big(\lam^2+ b^2(\bB_1^2+\bB_2^2)\big)+\vka_1b^2(\bS_1\bB_1+\bS_2\bB_2)\big]
=-(\ga_1+\lam)^2\lam\big(\lam^2+b^2(\bB_1^2+\bB_2^2)\big)
$$
since $\bS_1\bB_1+\bS_2\bB_2 =0$ by \eqref{sec}.
Hence, 
$$
\lam_1=\lam_2=-\ga_1, \quad \lam_3=0, \quad \lam_{4,5}=\pm ib\sqrt{\bB_1^2+\bB_2^2}.
$$
 II. For $(\bM,\bS)\in Z^r_2$, the Jacobian reads
  \be\la{J2}
J=\frac 12\begin{pmatrix}
-\ga_1 &0&  \vka_1 & 0 & 0
\\ 
0& -\ga_1 &0&   \vka_1 &0 
\\
-b \bS_3 & 0 & 0& 0& 0
\\
0&-b \bS_3  & 0& 0& 0
\\
b\bS_1 &b\bS_2& 0  &0& 0
\end{pmatrix}.
\ee
Hence,
$$
\det(2J-\lam)=
\begin{bmatrix}
-\ga_1-\lam & 0 & \vka_1&0&0
\\
0&-\ga_1-\lam & 0&\vka_1&0
\\
-b \bS_3&0&-\lam&0& 0
\\
0&-b \bS_3&0&-\lam &0
\\
b\bS_1& b\bS_2&0&0&-\lam
\end{bmatrix}
=-\lam\begin{bmatrix}
-\ga_1-\lam & 0 & \vka_1&0
\\
0&-\ga_1-\lam & 0&\vka_1
\\
-b \bS_3&0&-\lam&0
\\
0&-b \bS_3&0&-\lam 
\end{bmatrix}
$$
$$
=-\lam\big[(\ga_1+\lam) \begin{bmatrix}
\ga_1+\lam & 0&\vka_1
\\
0&-\lam&0
\\
b \bS_3&0&-\lam 
\end{bmatrix}
+b \bS_3
\begin{bmatrix}
0 & \vka_1&0
\\
\ga_1+\lam & 0&\vka_1
\\
b \bS_3&0&-\lam 
\end{bmatrix}\big]
$$  
$$  
= -\lam\big[(\ga_1+\lam)(\lam^2(\ga_1+\lam)+\lam  \vka_1 b\bS_3) +b \bS_3(\vka_1^2b\bS_3+\lam\vka_1(\ga_1+\lam) \big] 
$$
$$
= -\lam\big[\lam^2(\ga_1+\lam)^2+2\lam(\ga_1+\lam) b\vka_1 \bS_3 +b^2\vka_1^2\bS_3^2 \big]
=-\lam\big[\lam^2 +\lam\ga_1+b\vka_1 \bS_3\big]^2.   
$$
Hence, we have the following roots:
\be\la{eig2}
\lam_{1,2}(\bS_3)=-\frac{\ga_1}2+ \frac{\sqrt{\ga_1^2-4b\vka_1\bS_3}}2,\qquad \lam_{3,4}(\bS_3)=-\frac{\ga_1}2- \frac{\sqrt{\ga_1^2-4b\vka_1\bS_3}}2,\qquad \lam_5=0.
\ee
By (\ref{P1}),  the 
spectrum of the linearised system (\ref{perNa}) at the 
stationary states 
 $(\bM,\bS)\in Z^r_2$
consists of 
$$
p\lam_{1,2}(\bS_3)=-\frac{\ga}2+ \frac{\sqrt{\ga^2-4p^2b\vka_1\bS_3}}2,\qquad p\lam_{3,4}(\bS_3)=-\frac{\ga}2- \frac{\sqrt{\ga^2-4p^2b\vka_1\bS_3}}2,\qquad \lam_5=0.
$$
By (\ref{P1}), we have
 $b\vka_1=\fr{2\om}{\hbar}>0$.
Hence, for $p\ne 0$,
\be\la{lam14}
\rRe \lam_k(\bS_3)<0\qquad{\rm for}\quad k=1,\dots,4
\quad{\rm iff}\quad \bS_3>0.
\ee

    \setcounter{equation}{0}
   \section{Attraction to stable submanifold}
 Let us assume that $r$  is sufficiently large, so that
 \be\la{bet}
cr>|\bA^e|.
%
 \ee
 Then $Z_2^r\ne\emptyset$ and  $|\al_r\bA^e|=|\bA^e|/(cr)<1$, and 
 we can rewrite (\ref{Sr2}) as
  \be\la{Sr22}
Z^r_2=\{\bZ^r_2(\bS_3)=(-\bA^e,-\al_r \bA^e_1,-\al_r\bA^e_2,\bS_3):\bS_3\in[-\beta_r,\beta_r]\},\qquad \beta_r:=\sqrt{1-\al_r^2|\bA^e|^2}> 0.
\ee
For $0\le a\le b\le1$ and $ d>0$ denote 
\be\la{Sr2e}
Z^r_+(a,b)=\{\bZ^r_2(\bS_3):\bS_3\in[a,b]\},\qquad \cU^r_ d(a,b)=\{\bZ^r_2(\bS_3)+\bN:
\bN\bot TZ^r_2(\bS_3),\,\,\, \bS_3\in[a,b],\,\,\,
|\bN|\le  d\}.
\ee
 By (\ref{lam14}), for $ s>0$ and $p\ne 0$, we have
 \be\la{max}
\max_{\ba{c}1\le k\le 4\\ \bS_3\in [ s,\beta_r]\ea}\rRe \lam_k(\bS_3)<-\nu(s)<0.
 \ee

The
  tubular domains $\cU^r_ d(a,b)$
  with $0< a\le b<1$
   are attracted to $Z^r_+$ in the following sense.
 \bl\la{lat}
 Let   (\ref{bet}) hold,
 $ s\in (0,\beta_r/4)$ and $(\bM( t),\bS( t))$ be the solution to the
 averaged equations (\ref{perNa}) with the initial state
 $(\bM(0),\bS(0))\in D^r_d:=\cU^r_ d(2 s,\beta_r-2 s)$.
  Then for sufficiently
 small $ d,\mu>0$,
 \be\la{Drd2}
\max_{ t\in[0,p^{-1}]} \dist
 ((\bM( t),\bS( t)),Z^r_+)=\cO(d).
 \ee
 \el
 \bpr 
  The Jacobian (\ref{J2})  admits 
  a basis of eigenvectors and generalised eigenvectors
  $v_k(\bS_3)$ with $\bS_3\in (0,1]$ and $k=1,\dots,5$ 
 corresponding
 to eigenvalues  (\ref{eig2}).  
 
 The set $Z^r_+$ is one-dimensional manifold of stationary states of the vector field 
$(\ov f_r,\ov g_r)$. Hence, the tangent vectors to $Z^r_+$ at any point
belong to the kernel $\Ker J$. So, the tangent vectors correspond to 
the eigenvalue $\lam_5=0$. Therefore, the eigenvectors corresponding to the stable 
eigenvalues  (\ref{lam14}), are transversal to  $Z^r_+$.
 The vectors can be chosen piece-wise continuous in $\bS_3\in (0,\beta_r]$ since the multiplicity 
 of the eigenvalues (\ref{eig2}) is constant.
 Define the map 
  $\R^4\times (0,\beta_r]\to\R^5$
  by
  \be\la{dex}
  (x_1,x_2,x_3,x_4,x_5)\mapsto(\bM,\bS)=\bZ^r_2(x_5)+\sum_1^4 x_kv_k(x_5),
  \qquad x_5\in(0,\beta_r].
  \ee
  Let us check that the map is nondegenerate at every point
  $(0,0,0,0,x_5)$ with $x_5\in(0,\beta_r]$ which corresponds to $(\bM,\bS)=\bZ^r_2(x_5)\in Z^r_+$.
  Indeed, 
   by (\ref{Sr22}), the last column of the Jacobian $\pa(\bM,\bS)/\pa x$ is the vector
$ {\pa (\bM,\bS)}/{\pa x_5}=(0,0,0,0,1)=v_5(x_5)$. On the other hand for $k\le 4$, we have 
$ {\pa (\bM,\bS)}/{\pa x_k}=v_k(x_5)$. Hence, the columns of the Jacobian are linearly
independent. Hence, the coordinates $x(\bM,\bS)$ can be choosen  piece-wise
continuous in a neighborhood of
$Z^r_2( s,\beta_r- s)$.
By (\ref{dex}) and (\ref{Sr2}), we have
 $x_5=\bS_3$ for $x_1=x_2=x_3=x_4=0$.
 
 Let us  fix an $ s\in(0,\beta_r/4)$  and sufficiently
 small $ d>0$ such that  the coordinates $x$ are well defined in 
 $\cU^r_ d( s,\beta_r- s)$.
  Then
 the distance $d=d((\bM,\bS),Z^r_2)$
  for $(\bM,\bS)\in \cU^r_ d( s,\beta_r- s)$ 
  is equivalent to  $\bbd=\bbd(x)=\sqrt{\sum_1^4 x_k^2}$.
  For $\de>0$ let us    denote 
   $$
   \cT^r_\de( s)=\{x\in\R^5: x_5\in[ s,\beta_r- s],\bbd\le  \de\}. 
  $$
  Then 
   the coordinates $x$ are well defined in 
 $\cT^r_\de( s)$, and 
     the equations (\ref{perNa}) become
  \be\la{perNab}
  \left\{\ba{rclrcl}
 \dot x_1&=&p[\lam_1(x)x_1+q_1(x)],&
   \dot x_2&=&p[\lam_2(x)x_2+\ve_1(x)x_1+q_2(x)]
  \\\\
  \dot x_3&=&p[\lam_3(x)x_3+q_3(x)],& \dot x_4&=&p[\lam_4(x)x_4+\ve_2(x)x_3+q_4(x)],
 \qquad \dot x_5=pq_5(x),
  \ea\right|,
  \ee
  where 
  $q_k(x)=\cO(\bbd^2)$ and $|\ve_k(x)|$ with $k=1,2$ are sufficiently small.
     By (\ref{max}), the equations (\ref{perNab}) imply that for solutions 
     $x( t)\in \cT^r_\de( s)$ and sufficiently small $ \de>0$,
       \be\la{ddd}
 \pa_ t \bbd^2(x( t))\le 0.
 \ee
 Therefore, for any solution $x( t)$ with initial state $x(0)\in \cT^r_\de( s)$,
 \be\la{nnd}
 \bbd^2(x( t))\le  \de^2 \quad{\rm until}\quad x_5( t)\in [ s,\beta_r- s].
 \ee
 This inequality
  allows us to estimate the exit time $ t_*$ from  $\cT^r_ \de( s)$ of 
 solutions with initial states $x(0)\in \cT^r_ \de(2 s)$.
 Namely, 
  in the case $ t_*<\infty$
 and sufficiently small $ \de>0$, the coordinate $x_5( t)$
 must pass either the segment $[ s,2 s]$ or  $[\beta_r-2 s,\beta_r- s]$.
 In both cases the last equation of (\ref{perNab}) 
 together with
  (\ref{nnd}) imply that
 \be\la{exi}
  s/2\le \int_0^{ t_*} |\dot x_5( t)|d t
 \le  pC\int_0^{ t_*} \bbd^2(x( t))d t
 \le pC \de^2\int_0^{ t_*}d t=pC \de^2t_*.
  \ee
  Thus, 
 $$ t_*>  p^{-1}\fr{s}{2C\de^2}.
 $$
 Finally, 
 for sufficiently  small $\de>0$ we can choose $s=2C\de^2$. Then
 (\ref{nnd}) implies (\ref{Drd2}) for sufficiently  small $d>0$.
  \epr

 \br\rm
 By (\ref{max}),
   equations  (\ref{perNab}) with $k\le 4$
 imply an exponential approach to the manifold  $Z^r_+$. Hence, 
 the last equation with $k=5$, gives the exponential decay of $|\dot x_5|$.
 This is why the trajectory remains for a long time in a small neighborhood of 
 $Z^r_+$.
 \er

        \setcounter{equation}{0}
  \section{Single-frequency asymptotics}\la{ssin}
  
   In this section, we prove 
    the asymptotics (\ref{intriad}),  (\ref{intrias})  in the resonance case $\Om=\om$.
 By the relations (\ref{intr}),  the asymptotics are equivalent to the
  corresponding asymptotics of solutions $(\aM(t),\aS(t))$ to the interaction
  dynamics (\ref{perN}) which we prove
     applying the results of the  averaging theory \ci{SVM2007}. 
   To justify the application, we are going to check suitable properties of the 
   system (\ref{perN}).

   \subsection{KBM vector field}
 Let us denote $v_r=(f_r, g_r)$ 
 the vector field (\ref{fg})
   of the system  (\ref{perN}).
 It is easy to check that in  the case of  {\it almost periodic pumping} $A^e(t)$ and each
 $r\!\ne\!0$, 
  for any bounded   region $\cD\subset\aX$, we have
  \be\la{KBM}
  \sup_{(\bM,\bS)\in \cD}\fr1{T}
 \Big| \int_0^T[v_r(\bM,\bS, t)-\ov v_r(\bM,\bS)]d t\Big|
  \to 0,\qquad T\to \infty.
  \ee  
  Moreover, 
  $v_r$  is the Lipschitz continuous 
 vector  field in any bounded region $\cD\subset\aX$. 
 Hence, $v_r$
   is a
KBM (Krylov--Bogolyubov--Mitropolsky) vector field  in any bounded region $\cD$
according to \ci[Definition 4.2.4]{SVM2007}.

Furthermore,    for the {\it quasiperiodic pumping} (\ref{qpp}),
the formulas
 (\ref{fg}) implies that
    \be\la{KBM2}
  \de_\cD(p):=p
  \sup_{(\bM,\bS)\in \cD}
  \,\,
  \sup_{T\in[0,p^{-1}]}
 \Big| \int_0^T[v_r(\bM,\bS, t)\!-\!\ov v_r(\bM,\bS)]d t\Big|
=\cO(p),\,\,\, p\to 0,
  \ee    
  where $\de_\cD(p)$
  is  
the corresponding {\it order function} defined in \ci[Lemma 4.6.4]{SVM2007}.

  \br\la{rD}\rm
   For more general almost periodic pumping $A^e(t)$, the order function
  can be different \ci[Section 4.6]{SVM2007}.
  \er

  \subsection{The  asymptotics in the interaction picture}
Here we prove Theorem \ref{t1}.
\smallskip\\
i)   The asymptotics (\ref{intriad}) is equivalent to 
 \be\la{intriadeq}
 \max_{t\in [0,p^{-1}]}\Big[|\aM(t)-\bM|+|\aS(t)-\bS|\Big]=\cO({p^{1/2}}),
 \qquad { p\to 0},\quad{p}/\ga=r,\quad\qquad\qquad\qquad\qquad
   \ee 
  The initial state
  $(\aM(0),\aS(0))=(M(0),S(0))$ by  (\ref{intr}).
  Hence,
  $(\aM(0),\aS(0))=(\bM,\bS)$. However,  $(\bM,\bS)$ is a stationary state for the averaged system  (\ref{aver}), or equivalently, to (\ref{aveq}). Therefore,  Theorem 4.3.6 of \ci{SVM2007} implies
  (\ref{intriadeq}).
  Indeed,  both conditions 1 and 2 of the theorem hold in our case since
  \smallskip\\
  a) $v_r$  is a KBM-vector
  field in any bounded region $\cD\subset\aX $ with the order function (\ref{KBM2});
  \smallskip\\
  b)  solutions $(\aM( t),\aS( t))$ of the system (\ref{systim}) with ${p}/\ga=r$ and small $p>0$
  are uniformly  bounded  by the 
  a priori estimate (\ref{apr2}). 
\smallskip\\  
   Hence, (\ref{intriadeq}) is proved.
  \smallskip\\
ii) Similarly, the asymptotics (\ref{intrias}) is equivalent to 
 \be\la{intriaseq}
 \max_{t\in [0,p^{-1}]}\Big[ |\aM(t)-\bM_*|+|\aS(t)- \bS_*(t)|\Big]
 =\cO[{p^{1/2}}+ d],
 \qquad p\to 0,\quad{p}/\ga=r.\qquad
   \ee 
   Define the domain  $D^r_ d$ as in Lemma \ref{lat} and denote by
   $(\bM(t),\bS(t))$  the solution to the averaged system (\ref{aveq})
   with the initial value  $(\bM(0),\bS(0))=X(0)\in D^r_d$.
    Then  
 the same Theorem 4.3.6 of \ci{SVM2007} implies that
    \be\la{Bt51}
 \max_{ t\in[0,p^{-1}]} | (\aM( t),\aS( t))-(\bM( t),\bS( t))|=\cO(p^{1/2}),
 \quad p\to 0,\quad{p}/\ga=r.
  \ee 
  On the other hand,   
    Lemma \ref{lat} implies that 
   there exist  $(\bM_*(t),\bS_*(t))\in Z^r_+$ such that
   \be\la{intriaseq2}
\max_{t\in [0,p^{-1}]}\Big[ |\bM(t)-\bM_*(t)|+|\bS(t)- \bS_*(t)|\Big]
 =\cO(d).
   \ee 
  Now (\ref{intriaseq})  follows since $M_*(t)=-\bA^e$ 
  by (\ref{Sr2}).

  \appendix
  
  \protect\renewcommand{\thetheorem}{\Alph{section}.\arabic{theorem}}

\setcounter{equation}{0}
\section{The a priori bounds and well-posedness}\la{swp}
In this section, we 
prove the a priori bounds (\ref{apri})
 assuming that $A^e(t)\in C[0,\infty)$.
The bounds
imply
 the well-posedness of the MBE in
the phase space $\aX=\C\times \R^3$. 
 The density matrices (\ref{roS}) all are bounded by the conservation (\ref{Sc}) .
   Hence, it remains to prove 
the a priori estimates for the Maxwell amplitudes $(A(t),B(t))$.   
The following lemma is proved in 
\ci{KK2025sin}.

\bl
There exists a Lyapunov function
$V(A,B)$ such that 
\be\la{aprL}
a_1[A^2+B^2]\le V(A,B)\le a_2[A^2+B^2]\quad {\rm where}\quad a_1,a_2>0,
\ee
 and for solutions to (\ref{HMB2}), the function $V(t)=V(A(t),B(t))$ satisfies
 the inequality 
\be\la{Vder3}
\dot V(t)\le -\ga b
V(t)
+d\fr{p^2}\ga,\qquad t>0;
\quad b,d>0.
\ee
 \el

\bc
Solving the inequality (\ref{Vder3}),
we obtain:
\be\la{Vdec}
V(t)\le 
V(0)+ \fr db r^2,\qquad t\ge 0,\quad r=p/\ga.
\ee
Hence,
for solutions to (\ref{HMB2}) with $p/\ga=r$,
 the following bounds hold:
 \be\la{apr2}
A^2(t)+B^2(t)\le D_r(A^2(0)+B^2(0)),\qquad t\ge 0.
\ee
Now (\ref{apri}) follows.
\ec

\setcounter{equation}{0}
\section{The Maxwell--Bloch equations for pure and mixed states}\la{sMB}

There are various versions of the Maxwell--Bloch equations, see for instance  \ci{L1964} and \ci{A1985,AB1965}.
In this section we recall introduction of the MBE
for pure states and also for mixed states 
described by density matrix.

\subsection{The Maxwell--Bloch equations for pure states}

 In \ci{K2024}, the 
 MBE  for pure states were obtained as the Galerkin approximation of the
 damped driven
Maxwell--Schr\"odinger system. 
The approximation
consists of a single-mode Maxwell field coupled to
  two-level molecule in a bounded cavity $V\subset\R^3$:
 \be\la{solMB2}
\bA(x,t)=A(t)\bX(x),\quad
\psi(x,t)=C_{1}(t)\psi_{1}(x)+C_{2}(t)\psi_{2}(x),\qquad x\in V.
\ee
 Here 
 $\bA(x,t)$ denotes the vector potential of the Maxwell field, and
 $\bX(x)$ is a normalised eigenfunction of the Laplace
 operator in $V$ under suitable boundary value conditions
 with an eigenvalue  $-\Om^2/c^2$. By
$\psi_{l}$ we denote  some 
normalised eigenfunctions of the 
Schr\"odinger operator $\bH:=-\fr{\hbar^2}{2\cm} \Delta  + e\Phi (x)$ with the corresponding 
eigenvalues $\hbar\om_1<\hbar\om_2$, where $\Phi (x)$
is the molecular (ion's) potential.
The MBE  read as the Hamiltonian system with a dissipation and an external source:
\be\la{HMB} 
\fr1{c^2}\dot A(t)=\pa_B H,\quad \fr1{c^2}\dot B(t) =-\pa_A H-\fr\ga{c^2} B;
\qquad i\hbar\dot C_{l}(t)=\pa_{\ov C_{l}}H,\quad l=1,2.
\ee
Here the Hamiltonian is defined as
$$
H(A,B,C_1,C_2,t)=\cH(A\bX,
B\bX, C_{1}\vp_{1}
+C_{2}\vp_{2},t),
$$
where
$\cH$ is the Hamiltonian of the coupled
Maxwell--Schr\"odinger equations 
with pumping.
Neglecting the spin and scalar potential (which can be easily added),
the Hamiltonian $H$,  in the traditional {\it dipole approximation},
reads as \ci[(A.5)]{K2024}:
\be\la{Hc33}
H(A,B,C,t)=\fr1{2c^2}[B^2
+ \Om^2A^2]+\hbar\om_1 |C_{1}|^2+
\hbar\om_2|C_{2}|^2
  -\fr {2\vka}c
[A+A^e(t)]\,
 \rIm[\ov C_{1}C_{2}],\quad C=(C_1,C_2).
\ee
 Now the Hamilton equations (\ref{HMB}) become 
 \be\la{HMB21} 
\left\{\!\!\!\!\ba{rcl}
\dot A(t)\!\!\!&\!\!\!=\!\!\!&\!\!\!B(t),\quad
\dot B(t)=-\Om^2 A(t)\!-\!\ga  B(t)+cj(t)
\\\\
 i\hbar\dot C_{1}(t)\!\!\!&\!\!\!=\!\!\!&\!\!\!\hbar\om_1C_{1}(t)+ia(t)\,
  C_{2}(t),\,\,\,
i\hbar\dot C_{2}(t)=\hbar \om_2C_{2}(t)-ia(t)\,
  C_{1}(t)
 \ea\right|,\,\,\, j(t)={2}\vka\,\rIm[\ov C_{1}(t)C_{2}(t)].
\ee 
 The charge conservation 
 $|C_{1}(t)|^2+|C_{2}(t)|^2=\const$
 follows by differentiation from the last two equations of  (\ref{HMB21})
since the function $a(t)$ is real-valued. 
We  consider solutions with $\const=1$, cf. (\ref{Hro}):
\be\la{Hro2}
|C_{1}(t)|^2+|C_{2}(t)|^2=1,\qquad t>0.
\ee

\subsection{The von Neumann equation for mixed states}

The Schr\"odinger amplitudes $C_1(t)$ and $C_2(t)$ 
in (\ref{HMB2})
can be replaced by density matrix $\rho$
which is a nonnegative Hermitian  $2\times 2$-matrix (\ref{Hro}).
In particular, for 
 the wave function  (\ref{solMB2}), the corresponding density matrix reads as
 \be\la{dmw}
\rho(t)=
|C(t)\rangle\langle \ov C(t)|=
\left(\!\!\!\ba{cc}|C_1(t)|^2&C_1(t)\ov C_2(t)\\\\C_2(t)\ov C_1(t)&|C_2(t)|^2\ea\!\!\!\right),
\qquad C(t)=\left(\!\!\!\ba{c}C_1(t)\\C_2(t)\ea\!\!\!\right).
\ee
The matrix satisfies the trace condition in (\ref{Hro}) by the charge conservation (\ref{Hro2}). 
The last line of (\ref{HMB21})  in the vector form reads
 \be\la{dmw2}
i\hbar\dot C(t)=H(t)C(t),\qquad H(t)=
\left(\!\!\!\ba{cc}\hbar\om_1&ia(t)\\\\-ia(t)&\hbar\om_2\ea\!\!\!\right)\!\!.
\ee
Accordingly, 
the density matrix (\ref{dmw}) satisfies  von Neumann equation  \ci{SZ1997}:
 \be\la{dmw3}
\!\!\!\!\!\!\!i\hbar\dot\rho(t)=i\hbar
[\dot C(t)\otimes \ov C(t)+C(t)\otimes \dot{\ov C}(t)]
=
[H(t)C(t)]\otimes \ov C(t)-C(t)\otimes [H(t)\ov{C}(t)]
=[H(t),\rho(t)].
\ee
Now the Maxwell--Bloch system (\ref{HMB21}) is replaced by (\ref{HMB2}).

General density matrix describes
an ensemble of molecules with the pure states
$C(k,t)=(C_1(k,t),C_2(k,t))$ 
satisfying (\ref{Hro2}):
\be\la{denmar}
\rho(t)=\sum_k p_k C(k,t)\otimes\ov C(k,t),\qquad p_k\ge 0,\quad 
\sum_k p_k=1,
\ee
 where $p_k$ are the probabilities of the pure states.

\setcounter{equation}{0}
  \section{On possible treatment of the laser action}\la{aC}
Here we discuss possible treatment of the laser action 
  relying on the obtained results.
  \smallskip\\
  {\bf On the smallness of the parameters.}
Note that  the dipole moment $p$ and dissipation coefficient $\ga$ are very small 
for many types of lasers. In particular, 
the dissipation coefficient  for the Ruby laser is the electrical conduction of corundum
which is 
 $\ga\sim10^{-14}$ in the Heaviside--Lorentz 
 units
\ci{DH1992,S1986,S2012}.
For the dipole moment typically 
$|p|\sim 10^{-18}$ according to
\ci{MW}
that agrees with the classical dipole moment $ed/2$, where 
$d\sim 10^{-8}$cm is the molecular diameter, and
$e\sim 10^{-10}$ is the elementary 
charge in the same units. 
 \smallskip\\
{\bf On the laser threshold.}
The asymptotics (\ref{intriad}) hold for solutions 
to (\ref{HMB2})
with harmonic initial states $X(0)\in Z^r$, and the Lebesgue measure $|Z^r|=0$. On the other hand, (\ref{intrias}) hold for solutions with  initial states from an open
domain of attraction. Hence,
by (\ref{slo}),
the asymptotics (\ref{intrias}) appear with a ``nonzero probability" in contrast to 
 (\ref{intriad}). This fact, provisionally, clarifies the existence of a
 {\it laser threshold}
 to ignite the laser action: the   intensity of random pumping  must be sufficiently large to bring the solution to
  the domain
 of attraction, and then
  the solution is  captured in the domain with the single-frequency asymptotics.
 \smallskip\\
 {\bf On the laser amplification.}
The equations (\ref{HMB2}) describe  one molecule  coupled to
 the Maxwell field. 
The limiting amplitudes of the Maxwell field in the asymptotics 
(\ref{intriad}), (\ref{intrias})
   do not depend on non-resonance harmonics in the pumping (\ref{qpp})
  with the frequencies $\Om_k\ne\om$.
 This means that the dynamics  (\ref{HMB2})  acts as a filter, selecting only the resonant harmonics, that itself cannot explain the
  amplification of the Maxwell field in laser devices.
 The amplification could be explained by a large number of active  molecules, typically $N\sim 10^{20}$, 
 under the traditional assumption 
 that the molecules
 interact with the Maxwell field but do not interact with each other \ci{N1973}.
 In this case, the amplitude of the total Maxwell field is multiplied by $\sqrt{N}\sim 10^{10}$ by the Law of Large Numbers.
 
 This conclusion essentially depends on the fact that  
the phases of the amplitudes  $\bM^r$ and $\bM_*$
 in  asymptotics  (\ref{intriad})--(\ref{intrias})
for all molecules are  approximately uniformly distributed. 
The uniformity  takes place
 if it holds for  phases of  $\bA^e$.
 This follows from the asymptotics and 
Remark \ref{runif} since $e^{i\phi}e^{-i\Om t}\bM$ differs from $e^{-i\Om t}\bM$
by a shift of time.
 \smallskip\\
{\bf Self-induced transparency.} The harmonic states (\ref{Sr2})  means that the incident wave 
$\bA^e e^{-i\Om t}$ results in the outgoing wave $M(t)\approx -\bA^e e^{-i\Om t}$ with the same 
amplitude and frequency, but with the phase jump $\pi$. This phenomenon resembles 
 the self-induced transparency \ci{S1986}.

\end{document}